\documentclass[2025,preprint,12pt]{elsarticle}
\usepackage[margin=1in]{geometry}
\usepackage{amsmath,amssymb,amsthm,bm,mathtools}
\usepackage{graphicx}
\usepackage{xr-hyper}
\usepackage{verbatim}
\usepackage{xcolor}
\usepackage[colorlinks=true,citecolor=blue,linkcolor=blue,urlcolor=blue]{hyperref}

\usepackage{caption,subcaption}
\usepackage{float}
\usepackage{todonotes}

\newcommand{\R}{\mathbb{R}}

\biboptions{numbers,sort&compress}

\begin{document}

\begin{frontmatter}

\title{Singularities in Multi-Objective Optimization \\ and their Crossing during Continuation}

\author[1]{Arjun Manoj}
\author[2]{Michail E. Kavousanakis}
\author[3]{Shanqing Liu}
\author[1,4]{Ioannis G. Kevrekidis\corref{cor1}}
\ead{yannisk@jhu.edu}

\address[1]{Chemical \& Biomolecular Engineering, Johns Hopkins University, Baltimore, MD 21218, USA}
\address[2]{School of Chemical Engineering, National Technical University of Athens, Zografou, Athens 15780, Greece}
\address[3]{Division of Applied Mathematics, Brown University, Providence, RI 02906, USA}
\address[4]{Applied Mathematics and Statistics, Johns Hopkins University, Baltimore, MD 21218, USA}
\cortext[cor1]{Corresponding author.}

\begin{abstract}

Continuation methods help trace Pareto sets in multi-objective optimization but are inherently local: a single run traces a single connected branch, requiring multiple restarts to recover disconnected components of Pareto fronts. We show that, for unconstrained bi-objective problems under weighted-sum scalarization, these disconnects can be artifacts of singularities in the scalarization parameter, where the weight $\lambda$ diverges as the objective gradients become collinear. 
Recasting Pareto optimality as a nonlinear system, we apply pseudo-arclength continuation to follow the Pareto-critical set, and
show that suitable singular reparameterizations allow crossing these singularities in systematically, recovering disconnected branches in a single run. A coordinate-wise projective compactification further provides a unified framework for parameter and decision-space variables. We demonstrate the approach on the ZDT3 benchmark and modifications.

\end{abstract}

\begin{keyword}
Multi-objective optimization \sep Numerical continuation \sep Singularities \sep Pareto front \sep Compactification
\end{keyword}

\end{frontmatter}

\section{Introduction}\label{sec:intro}
Many problems in science and engineering involve optimizing multiple competing objectives, where improvement in one metric often degrades another. The set of trade-offs between objectives is captured by the \emph{Pareto set}, whose image in objective space is the \emph{Pareto front} \cite{miettinen1999nonlinear}. %deb2001}.
Scalarization methods combine multiple objectives into a single one, typically through weighted sums; while convenient, they recover only {\em convex regions} of the front and may miss non-convex trade-offs~\cite{das1997}. Alternatives, such as NBI~\cite{das1998}, 
and $\varepsilon$-constraint formulations~\cite{ehrgott2008epsilon} address nonconvex regions but introduce additional constraints and associated constraint-handling procedures.
%
%\arjun{More complete literature review with more citations}
%

Continuation methods~\cite{doedel1981auto} provide a geometric framework: building on first-order Pareto criticality conditions, they trace a manifold of Pareto-critical solutions through predictor--corrector schemes~\cite{hillermeier2012nonlinear}, sometimes combined with global search strategies or trajectory-based formulations~\cite{schutze2008hybridizing,martin2018pareto,schutze2025pareto,potschka2011,bolten2021}. 
Despite these advances, continuation-based methods remain inherently local: a single run traces only the \textit{connected component reachable from its initialization}, and may become ill-conditioned or fail when the scalarization parameterization becomes singular. As a result, complete recovery of a disconnected Pareto front typically requires multiple initializations and manual restarts.
%
%In addition, Potschka et al.~\cite{potschka2011} combined the normalized normal constraint (NNC) formulation with continuation-based integration schemes to trace Pareto-critical trajectories. 
%

A different perspective comes from nonlinear dynamics: in models exhibiting finite-time blow-up, apparent singularities can be \emph{traversed} through singular reparameterization or \textit{compactification}~\cite{kevrekidis2017infinity, kavousanakis2025goingflowsolvingsymmetrydriven}. In this Letter, we connect this idea to continuation-based MOO. 
%%YGK this is the place to cite our paper once on arXiv or submitted 
%
%
Using classical weighted scalarization for an unconstrained bi-objective problem, we show that the disconnects seen under weighted-sum continuation are induced by \textit{singularities in the parameterization itself}, rather than by the underlying function landscape, and that traversing them recovers the full disconnected Pareto-critical structure within a single continuation run. 
We demonstrate this on the ZDT3 benchmark and its modifications; Additionally we show that the construction generalizes from the scalarization parameter also to singularities in the decision variables via an alternative projective compactification. 
Throughout this work, we restrict attention to disconnections arising from singular continuation parameterizations, leaving intrinsically disconnected Pareto-critical manifolds for future investigation.

\section{Continuation of Pareto optimal solutions}\label{sec:method}
Consider the unconstrained MOO problem: $\min_{\bm{x}\in \mathcal{X}} \bm{f}(\bm x) = (f_1(\bm x),\dots,f_m(\bm x))$, where $\bm{x} \in \R^n$ denotes the decision vector, $\mathcal{X}$ is the feasible set, and $f_i:\R^n \to \R$ for $i=1,\dots,m$ are the multiple objective functions. A point $\bm{x}^\star$ is \emph{Pareto-optimal} if there exists no $\bm x \in \mathcal{X}$ such that $f_i(\bm x)\le f_i(\bm x^\star)$ for all $i$, with strict inequality for at least one objective $f_j$. The collection of all such points forms the \emph{Pareto set}, whose image in the objective space is referred to as the \emph{Pareto front}.
A common approach reformulates the problem through scalarization, applying first-order (KKT) optimality  conditions~\cite{hillermeier2012nonlinear, rakowska1993multi}. A point $\bm{x}^*$ is \emph{Pareto critical} if there exists a weight vector $\bm{\lambda} \in \mathbb{R}^m$, with $\lambda_i \ge 0$ and $\sum_{i=1}^k \lambda_i = 1$, such that $\sum_{i=1}^m \lambda_i \nabla f_i(\bm{x}^*) = \bm{0}$.
For the bi-objective case ($m=2$), weighted-sum scalarization reduces to: $g(\bm{x},\lambda) = \lambda f_1(\bm{x}) + (1-\lambda) f_2(\bm{x})$, $\lambda \in [0,1]$, which yields the stationarity condition
\begin{equation}
    \bm F(\bm x,\lambda) :=\nabla_{\bm{x}} g(\bm{x},\lambda) = \bm{0}, \qquad\text{for some } \lambda \in [0,1].
\end{equation}
This defines the Pareto-critical manifold
$\mathcal{S}=\{(\bm x,\lambda):\bm F(\bm x,\lambda)=\bm 0\}
$
in the augmented $(\bm x,\lambda)$-space. Every Pareto-optimal point is Pareto-critical and therefore lies on $\mathcal{S}$ for some $\lambda\in[0,1]$, although the converse need not hold. 
To trace $\mathcal{S}$, we use pseudo-arclength continuation (see~\ref{app:pac}). 
%Section~\ref{SI-app:pac} of the Supplementary Material).
While the classical approach follows the Pareto-critical manifold for $\lambda\in[0,1]$; we instead continue beyond this interval to investigate whether disconnected Pareto regions can be connected within a single run. As shown next, this continuation encounters singularities where $\lambda$ diverges.

\section{Results and Discussion}\label{sec:results}
\subsection{Appearance of $\lambda$-singularities}\label{sec:appearance}

We consider the standard ZDT3 benchmark problem in bi-objective form: 
\begin{equation}    
 f_1(\bm x)=x_1, \quad       
 f_2(\bm x)=\ell(\bm x)\,h(f_1,\ell),
\end{equation}
with $\ell(\bm x)=1+\tfrac{9}{d-1}\sum_{i=2}^d x_i$ and $h(f_1,\ell)=1-\sqrt{f_1/\ell}-(f_1/\ell)\sin(10\pi f_1)$. 
Restricting to the 2D case $(x,y)\in[0,1]^2$ and noting that $\ell$ is minimized in $y$ at $y=0$ (see 
\ref{app:zdt_2d}),
%Section~\ref{SI-app:zdt_2d} of the Supplementary Material),
the stationarity reduces to the 1D condition
\begin{equation}\label{eq:zdt3}
F(x,\lambda):=\partial_x g(x,0;\lambda)=\lambda f_1'(x) + (1-\lambda) f_2'(x)=0.
\end{equation}

\begin{figure}[t]
\centering
\includegraphics[width=\linewidth]{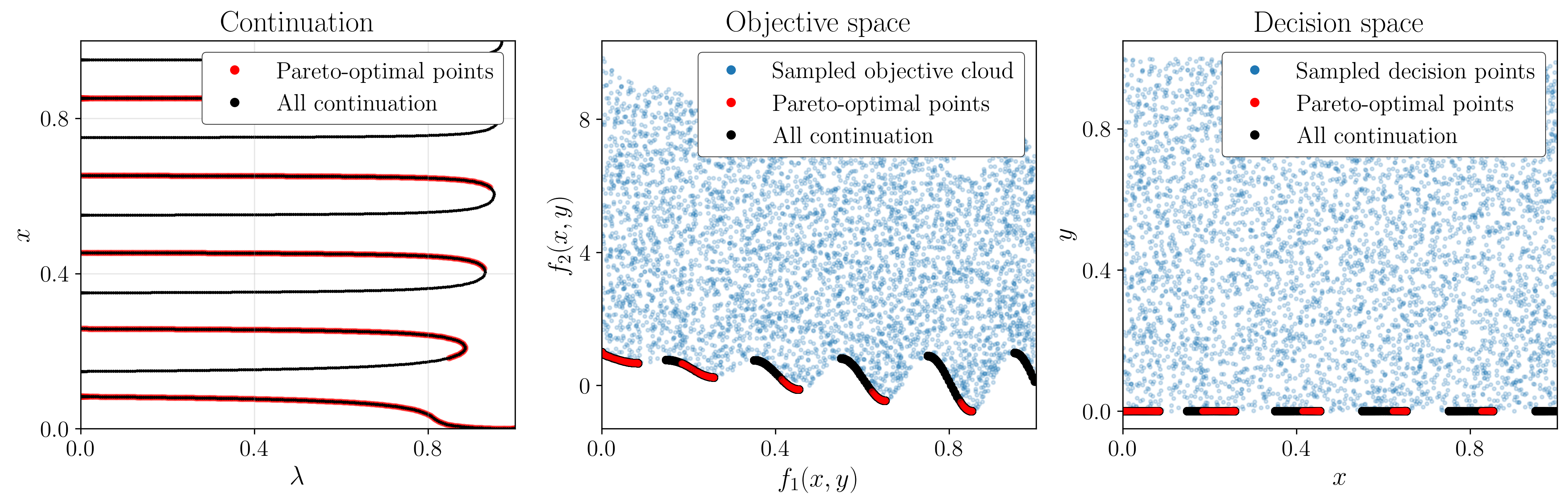}
\vspace{-0.5em}
\caption{Disconnected Pareto front of ZDT3. Purple points correspond to continuation branches obtained from multiple initializations; red are the Pareto-optimal points filtered by non-dominated sorting. (\textit{Left}) Pseudo-arclength continuation in $(\lambda,x)$ showing multiple branches of $g_x(x;\lambda)=0$. (\textit{Center}) Projection of the continuation branches into the objective space $(f_1(x, y),f_2(x, y))$. (\textit{Right}) Projection of the continuation branches into the decision space $(x,y)$.}
\label{fig:zdt32d}
\end{figure}

Tracing $F(x, \lambda)=0$ by pseudo-arclength continuation in $(x,\lambda)$ fails to recover certain boundaries of the feasible objective space, producing multiple disconnected branches (Fig.~\ref{fig:zdt32d}), which can only be obtained through manual restarts.
Rewriting Eq.~\ref{eq:zdt3} as: $f_1'({x}) = -\frac{1-\lambda}{\lambda}  f_2'({x})$, shows that the gradients must be collinear.
The gaps in Fig.~\ref{fig:zdt32d}, correspond to regions where \textit{increases} in $f_1$ accompany \textit{increases} in $f_2$; there the coefficient $\frac{1-\lambda}{\lambda}$ is negative, so $\lambda \notin[0,1]$ and the points cannot be Pareto-optimal. 

This raises a natural question: can continuation be extended beyond the globally Pareto-optimal (red) segments to recover the full continuation branch, including the dominated Pareto-critical (black) segments with $\lambda\in[0,1]$ and the missing connecting branches with $\lambda\notin[0,1]$—thereby connecting the disconnected Pareto regions—in a single run?
Solving $F(x, \lambda)=0$ for $\lambda$ gives $\lambda(x) = \frac{f_2'(x)}{f_2'(x)-f_1'(x)}$, which diverges as $f_2'(x)\to f_1'(x)$. These singularities in $\lambda$ are what disconnects the apparent branches.

\subsubsection{Compactification: $\lambda = \tan\theta$}\label{sec:tan}

\begin{figure}[H]
\centering
\includegraphics[width=\linewidth]{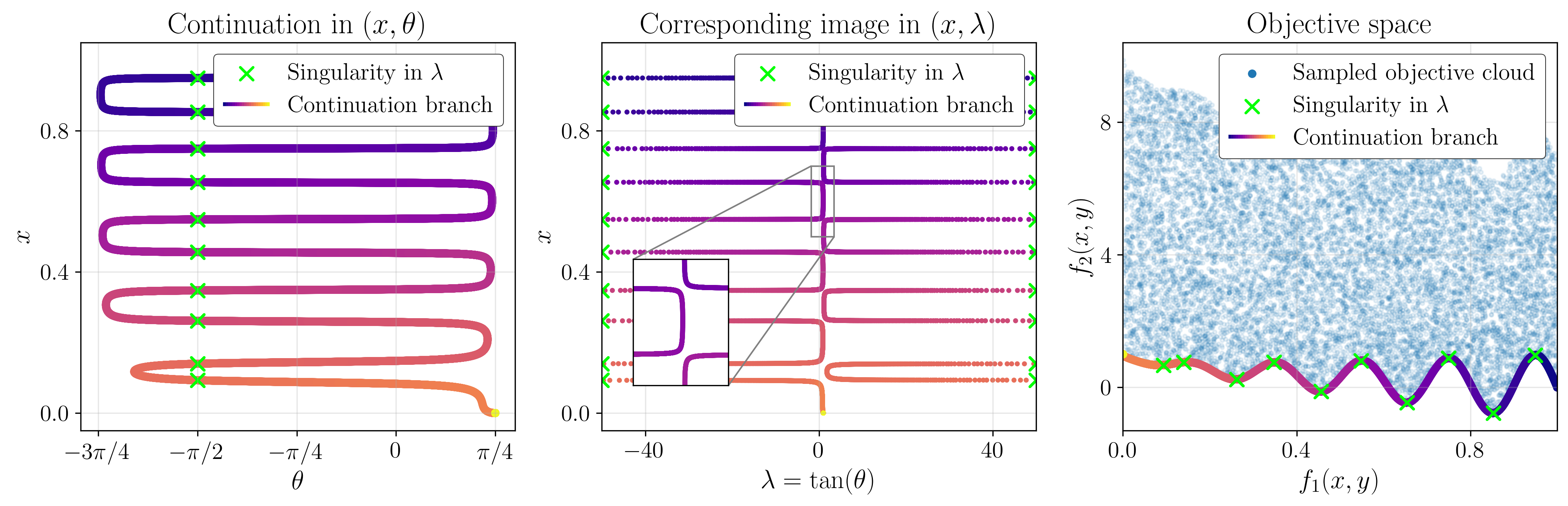}
\caption{Pseudo-arclength continuation under $\lambda=\tan\theta$ for ZDT3. \emph{Left:} smooth connected branch in $(\theta,x)$. \emph{Center:} the same branch in $(\lambda,x)$ appears fragmented. The inset highlights disconnections in $(x, \lambda)$. \emph{Right:} projection into objective space. Green markers indicate $|\lambda|$-singularities.}
\label{fig:theta}
\end{figure}

A natural fix is to allow $\lambda\in\R$ and compactify. Setting $\lambda=\tan\theta$ maps $\lambda\to\pm\infty$ to the finite point $\theta\to\pm\pi/2$. Implicit differentiation of $F=0$, combined with the stationarity condition $f_1'-f_2'=-f_2'/\lambda$, gives
\begin{equation}\label{eq:dxdlam}
\frac{dx}{d\lambda} = -\frac{-f_2'(x)}{\lambda^2 f_1''(x)+\lambda(1-\lambda)f_2''(x)} \;\xrightarrow{|\lambda|\to\infty}\; \frac{1}{c\lambda^2},
\quad c = f_1''(x^\star)-f_2''(x^\star),
\end{equation}
so $dx/d\lambda\to 0$. Under $\lambda=\tan\theta$, however, $dx/d\theta=(dx/d\lambda)\sec^2\theta$ with $\sec^2\theta=1+\lambda^2$, so
$\frac{dx}{d\theta} \xrightarrow{|\lambda|\to\infty} -\frac{1}{c}$, a finite limit. 
The $\infty/\infty$ indeterminacy is resolved by \textit{exact leading-order cancellation}. Crucially, $\theta$ is not restricted to $[-\pi/2,\pi/2]$: continuation lets $\theta$ evolve freely in $\R$, while $\tan\theta$ jumps repeatedly between $-\infty$ and $+\infty$ at each $\pi/2$-crossing. The full Pareto-critical manifold then appears as a single smooth curve in $(\theta,x)$ (Fig.~\ref{fig:theta} left panel); mapped back to $(\lambda,x)$, the same curve recovers the apparently disconnected branches of Fig.~\ref{fig:zdt32d}. The derivation extends to general $f_1$ with $f_1''\neq 0$ 
(\ref{app:tantheta}; see also Figs.~\ref{fig:zdt3x3}--\ref{fig:zdt3sinx}).
%(Section~\ref{SI-app:tantheta} of the Supplementary Material; see also Figs.~\ref{SI-fig:zdt3x3}--\ref{SI-fig:zdt3sinx}).

\subsubsection{Crossing infinity in $\lambda$ directly}\label{sec:inf}

\begin{figure}[H]
\centering
\includegraphics[width=\linewidth]{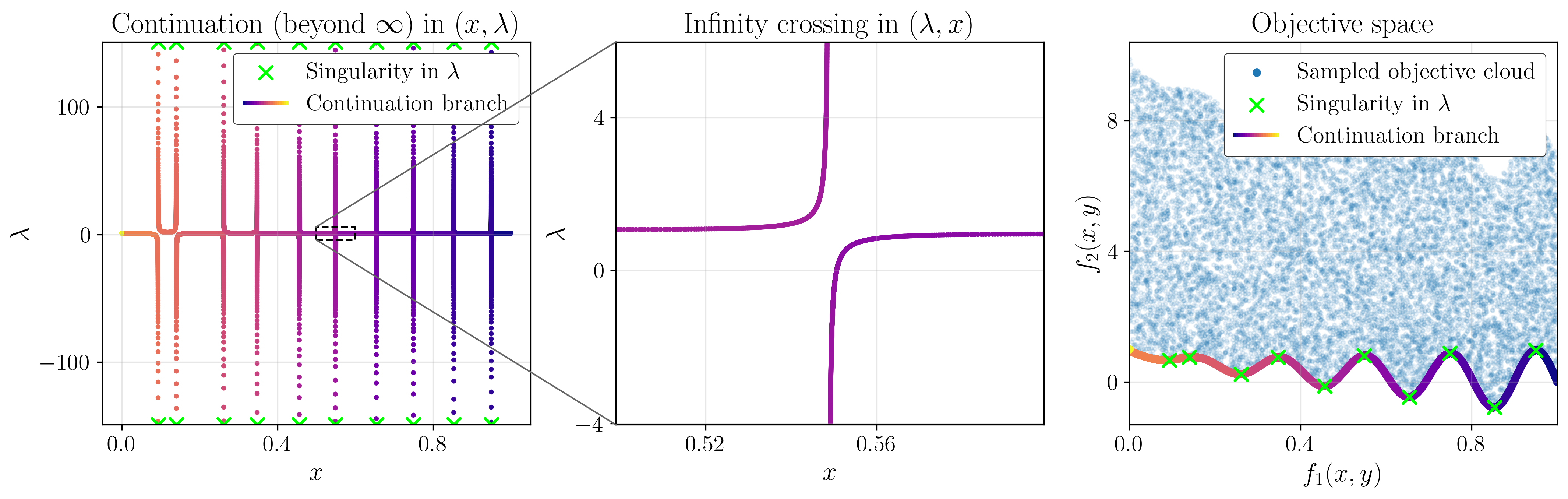}
\caption{Infinity crossing in $(\lambda,x)$ for ZDT3. \emph{Left:} branches recovered by reparameterization $v=1/(c\lambda)$. \emph{Center:}  magnified view of the singularity crossing obtained by integrating the regular ODE $dv/dx=-1$ through $v=0$. \emph{Right:} projection into objective space. Green markers indicate $\lambda$ singularities.}
\label{fig:inf}
\end{figure}

An alternative is to handle the singularity \emph{in situ}: we can cross $|\lambda|=\infty$ in place, following the construction of~\cite{kevrekidis2017infinity}. 
Near the singular point $x^\star$ (where $f_2'(x^\star)=f_1'(x^\star)$), Eq.~\eqref{eq:dxdlam} integrates to the asymptotic
\begin{equation}\label{eq:asym}
\lambda(x) \approx \frac{1}{c(x^\star-x)}, \qquad c=-\frac{f_2''(x^\star)-f_1''(x^\star)}{f_2'(x^\star)}.
\end{equation}
Once $|\lambda|>K\gg 1$ is detected, we record $N\ge3$ further continuation samples $(x_i,\lambda_i)$ and locally fit the approximation $x = x^\star - 1/(c\lambda)$ for $(x^\star,c)$ by least squares. Defining $v=1/(c\lambda)$ yields the regular ODE $dv/dx=-1$, which can be integrated smoothly through \(v=0\), corresponding to crossing the singularity at \(x=x^\star\). For the first solution point beyond the singularity, $x^{(0)}=x^\star+(x^\star-x_N)$, Eq.~\eqref{eq:asym} gives the clean reflection $\lambda^{(0)}=-\lambda_N$, a robust restart for pseudo-arclength on the next branch (Fig.~\ref{fig:inf}).

\subsection{Projective Compactification for Parameter  \underline{and} Decision Spaces}

While the singularities encountered here arise in the scalarization parameter $\lambda$, multi-objective continuation problems may also involve large or unbounded decision variables~\cite{lovison2011}, motivating a compactification framework that treats both on equal footing.
To illustrate this idea, we consider a modified two-dimensional ZDT3 problem that exhibits repeated $\lambda$-singularities together with large variations in a decision variable.
\begin{equation}\label{eq:modzdt3}
f_1(x,y)=x+y, \quad
f_2(x,y)=\ell\left(1-\sqrt{\tfrac{x}{\ell}}-\tfrac{x}{\ell}\sin(10\pi x)\right), \quad \ell=1+9y,
\end{equation}
on the extended domain $x\in[0,1.5]$, $y\in[-\tfrac19,\infty)$.
Along the Pareto-critical branch the oscillatory term drives the objective gradients into collinearity at several points, so $\lambda$ diverges as in Section~\ref{sec:appearance}.
The modified objective and enlarged feasible domain also generate regions where the continuation branch extends to very large values of the decision variable $y$. These features highlight the need for a compactification strategy capable of handling both parameter singularities and, more generally, large or unbounded excursions in decision space.

\begin{figure}[H]
    \centering
    \includegraphics[width=\linewidth]{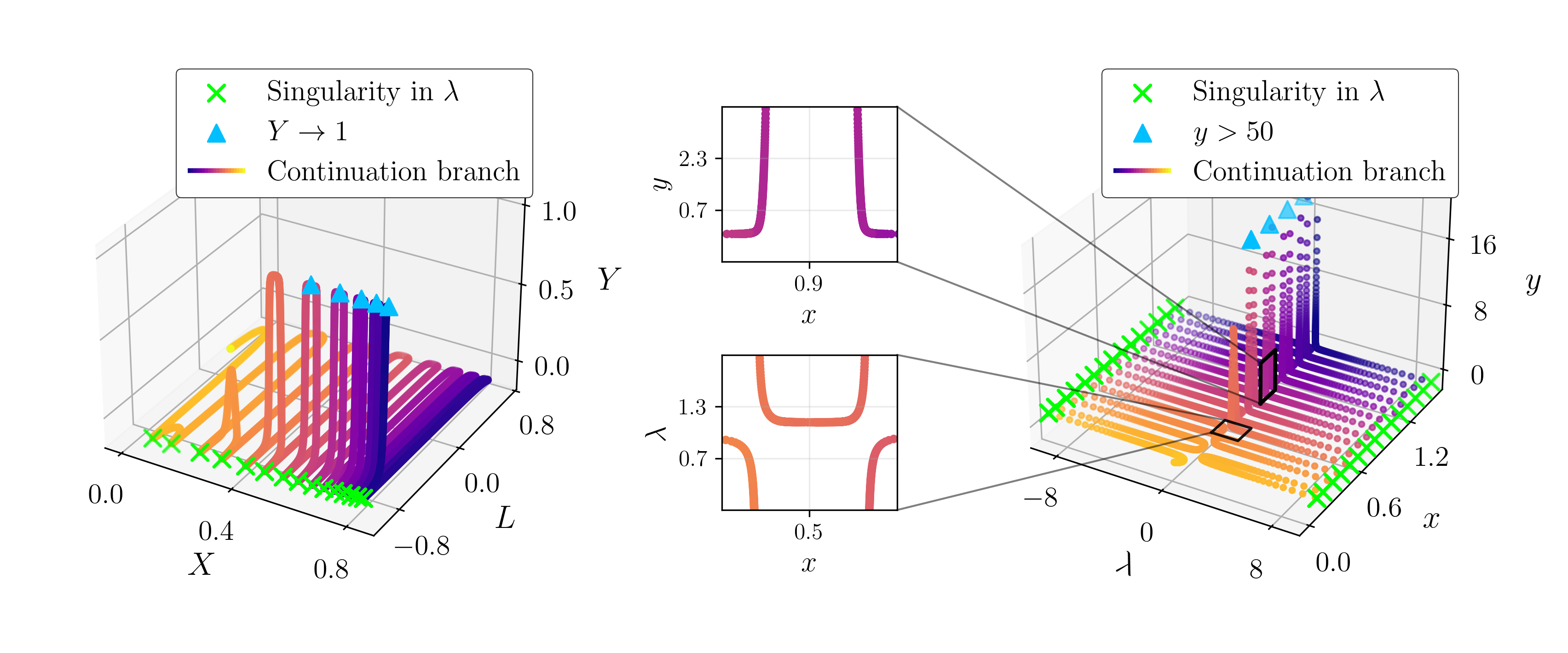}
    \includegraphics[width=\linewidth]{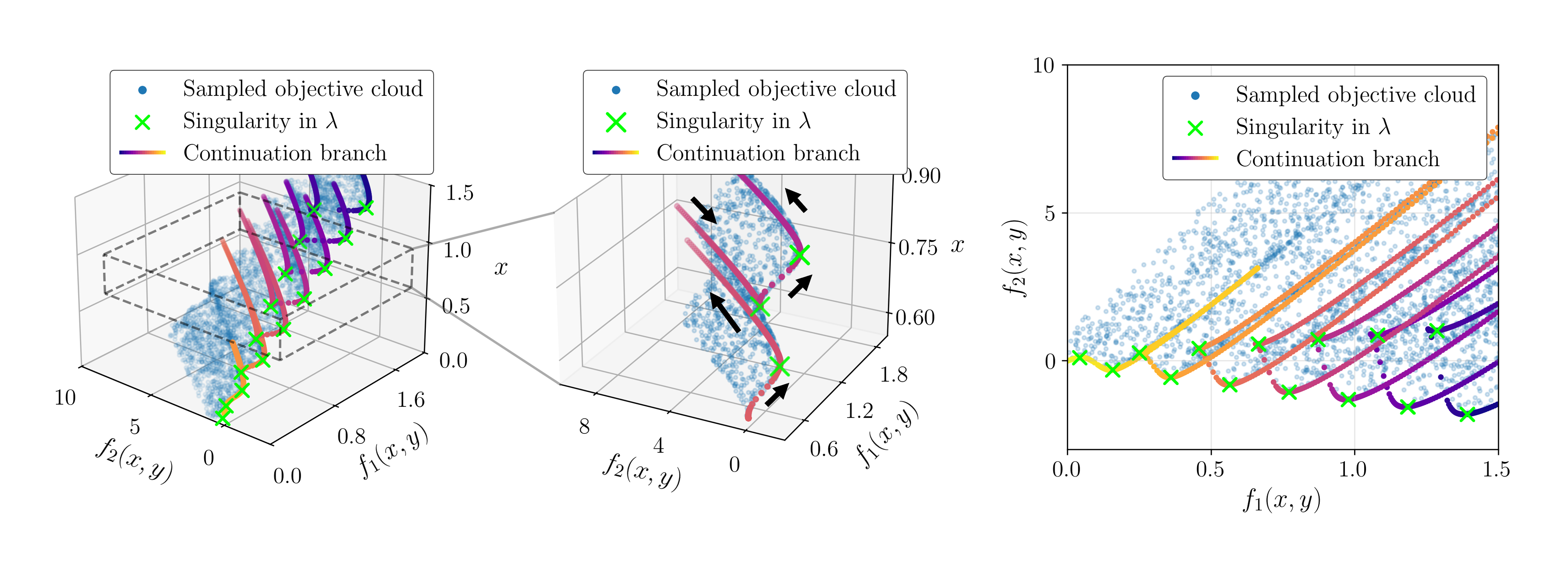}
    \caption{
    Pareto-critical continuation for the modified two-dimensional ZDT3 problem using coordinate-wise projective compactification. The continuation branch is shown in compactified coordinates $(X,Y,L)$ (top left) and mapped back to the original variables $(x,y,\lambda)$ (top right). The insets highlight regions of rapidly increasing $y$ (top) and the $\lambda$-singularities (bottom). While $y$ remains finite in this example, its successive peaks drive the compactified coordinate $Y$ progressively closer to the projective boundary $Y=1$. 
    The distinct $|\lambda|$-singularities are also traversed, with distinct infinities remaining separated under the coordinate-wise compactification. 
    The lower row visualizes the image of the continuation branch in the augmented objective space $(f_1,f_2,x)$(left), where the additional coordinate x separates branches that appear to overlap in the two-dimensional objective-space projection. The magnified, rotated view (center) shows how the continuation trajectory leaves the plotted region through the lower part of each wave and re-enters through the upper part (black arrows). These excursions remain finite, although their magnitude exceeds the plotted range. The corresponding projection into the objective space $(f_1,f_2)$ is shown on the right.
    }
    \label{fig:zdt3_pcompact}
\end{figure}

A single shared projective denominator would send every infinite limit to the same boundary, erasing the distinction between parameter and state singularities.
To keep them separate, we assign each variable its own projective chart, $X=\rho_x x$, $Y=\rho_y y$, $L=\rho_\lambda \lambda$, with independent normalizations: $X^2+\rho_x^2=1$, $Y^2+\rho_y^2=1$, and $L^2+\rho_\lambda^2=1$.
The compactified stationarity system remains regular at the boundary points $(\rho_x,X)=(0,\pm1)$, $(\rho_y,Y)=(0,\pm1)$, and $(\rho_\lambda,L)=(0,\pm1)$, now representing \emph{distinct} infinities of $x$, $y$, and $\lambda$. Pseudo-arclength continuation thus passes through both parameter \underline{and} state infinities of Eq.~\eqref{eq:modzdt3} in a single run, while preserving information about which variable diverges and where.
Figure~\ref{fig:zdt3_pcompact} shows the resulting branch in compactified coordinates, mapped back to the original variables, and projected into objective space, recovering the full Pareto-critical structure. Singularities confined to decision space alone, together with the single-scale construction, are illustrated for an additional bi-objective problem \cite{lovison2011} in \ref{app:decision_proj}. 
%Supplementary Material (Section~\ref{SI-app:decision_proj}).

\section{Conclusions}\label{sec:conclusion}
We have shown that the disconnected Pareto fronts exhibited by ZDT3 are not intrinsic to the optimization problem, but rather artifacts of the weighted-sum scalarization parameter: as the objective gradients become collinear, the weight $\lambda$ diverges and standard pseudo-arclength continuation fails. These are \textit{removable coordinate singularities} rather than folds or genuine breaks in the Pareto-critical manifold. Two simple reparameterizations cross them: angular compactification $\lambda=\tan\theta$ maps $\lambda \to \pm\infty$ to the finite point $\theta \to \pm\pi/2$, with the $\infty/\infty$ indeterminacy resolved by leading-order cancellation; continuation-through-infinity $v=1/(c\lambda)$ crosses the singularity in place via the regular ODE $dv/dx=-1$, yielding a clean reflection that restarts pseudo-arclength continuation on the next branch. Both depend only on second-derivative data at the singular point and reconstruct the full disconnected ZDT3 front in a single run. The same idea extends to decision-space singularities via coordinate-wise projective compactification. These compactifications do not eliminate genuine topological disconnections that arise from the underlying landscape; extending singularity-traversal ideas to such cases remains an interesting direction for future work.

\section*{Acknowledgements}
AM and IGK acknowledge support from the U.S.\ Department of Energy, Office of Science, Office of Advanced Scientific Computing Research (Award DE-SC0024162) and the National Science Foundation (Grant 2436738).

\bibliographystyle{elsarticle-num}
%\bibliographystyle{unsrt}
%\bibliography{biblio}

\appendix
\section{Pseudo-arclength continuation}\label{app:pac}
For $\bm F(\bm x,p)=\bm 0$ with $\bm F:\R^{n+1}\to\R^n$, natural-parameter continuation fails at folds, where the branch is not single-valued in $p$. Pseudo-arclength continuation \cite{keller1987lectures, kelley1995iterative, doedel1981auto} parameterizes the branch by an arclength variable $s$, predicts via $\bm x^{(0)}=\bm x_k+\Delta s\,\dot{\bm x}_k$, $p^{(0)}=p_k+\Delta s\,\dot p_k$, and corrects by Newton iteration on the augmented system
\begin{equation}
\bm F(\bm x,p)=\bm 0,\qquad (\bm x-\bm x_k)^\top\dot{\bm x}_k+(p-p_k)\dot p_k-\Delta s=0.
\end{equation}
The arclength constraint fixes $(\bm x,p)$ jointly on a hyperplane perpendicular to the tangent and traverses folds with $dp/ds=0$.

\begin{figure}[H]
    \centering
    \includegraphics[width=\linewidth]{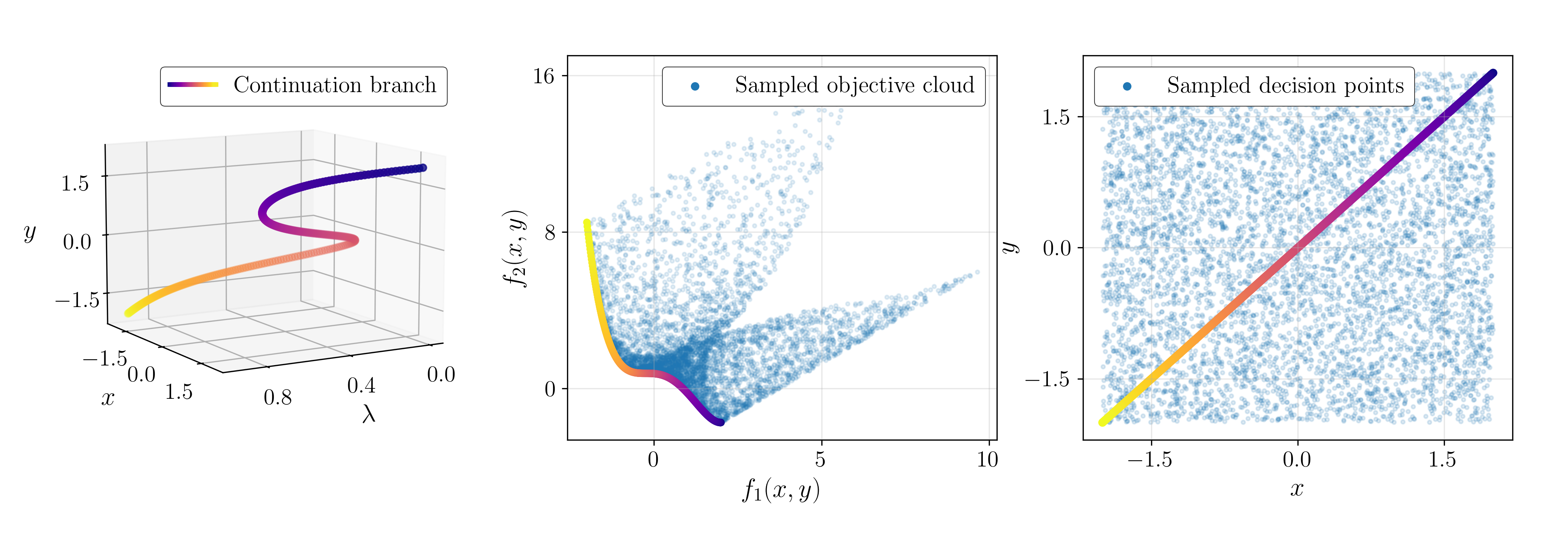}
    \caption{
    Pseudo-arclength continuation for a bi-objective problem with a non-convex Pareto front. (\textit{Left}) Continuation branch in the augmented $(x,y,\lambda)$ space. (\textit{Center}) Projection into the objective space $(f_1,f_2)$, illustrating recovery of the non-convex Pareto front. (\textit{Right}) Projection into the decision space $(x,y)$.
    }
    \label{fig:psarc}
\end{figure}

\section{Higher-dimensional extensions of ZDT3 continuation}\label{app:higher_dim}
The continuation framework extends beyond the one-dimensional, boundary-reduced ZDT3 example to higher-dimensional stationary manifolds arising from unconstrained ZDT3 formulations.

\subsection{Unconstrained two-dimensional ZDT3}\label{app:zdt_2d}
\begin{figure}[H]
    \centering
    \includegraphics[width=\linewidth]{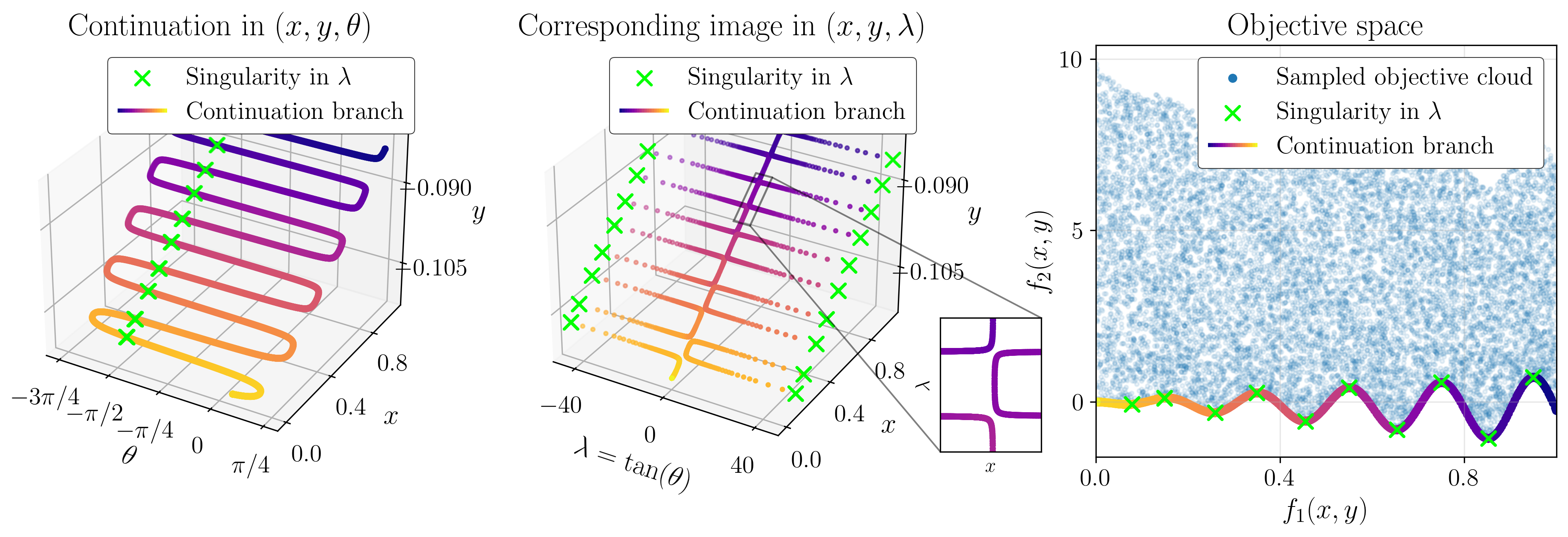}
    \caption{Unconstrained two-dimensional ZDT3 (no boundary reduction). Continuation in $(x,y,\theta)$ with $\lambda=\tan\theta$ remains smooth across $|\lambda|\to\infty$, while the image in $(x,y,\lambda)$ fragments. The inset highlights disconnections in $(x, \lambda)$.}
    \label{fig:zdt32dtheta}
\end{figure}

Consider the two-dimensional ZDT3 problem
\begin{equation}
f_1(x,y)=x, \qquad
f_2(x,y)=\ell(x,y)\left(1-\sqrt{\tfrac{x}{\ell(x,y)}}-\tfrac{x}{\ell(x,y)}\sin(10\pi x)\right),
\end{equation}
with $\ell(x,y)=1+9y$. In the standard constrained formulation the feasible region is $x\in[0,1]$, $y\in[0,1]$. The stationarity condition in $y$,
\begin{equation}
\frac{\partial f_2}{\partial y}=9-\frac{9}{2}\sqrt{\frac{x}{1+9y}}=0,
\end{equation}
admits no feasible interior roots, so all Pareto-optimal stationary points collapse onto the boundary $y=0$ and the continuation problem reduces to a 1D curve in $(x,\lambda)$. Relaxing the lower bound on $y$ makes interior stationary points feasible; solving the condition above gives
\begin{equation}
y=\frac{x/4-1}{9},
\end{equation}
which lies below zero for $x\in[0,1]$. Extending the domain to include these values, continuation proceeds in the full variable space $(x,y,\theta)$ with $\lambda=\tan\theta$, on the stationary manifold $\nabla_{x,y}[\lambda f_1+(1-\lambda)f_2]=0$. This produces a smooth 1D continuation curve embedded in 3D space that remains continuous across the singularities $|\lambda|\to\infty$, whereas its image in $(x,y,\lambda)$ shows apparent discontinuities and infinity crossings (Fig.~\ref{fig:zdt32dtheta}). The singularity is thus a coordinate artifact of the non-compact scalarization parameter, not an intrinsic feature of the problem.

\subsection{Modified two-dimensional ZDT3}\label{app:zdt_2d_mod}

\begin{figure}[H]
    \centering
    \includegraphics[width=\linewidth]{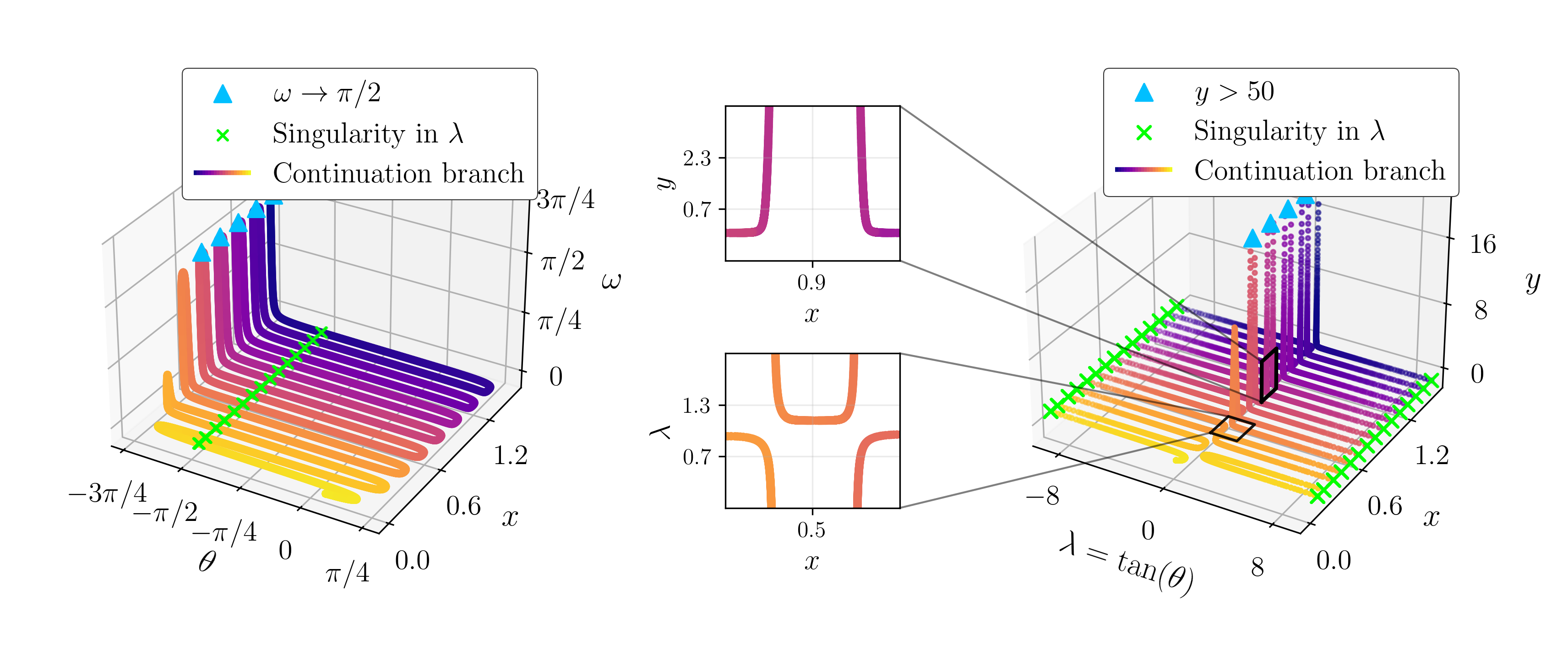}
    \includegraphics[width=\linewidth]{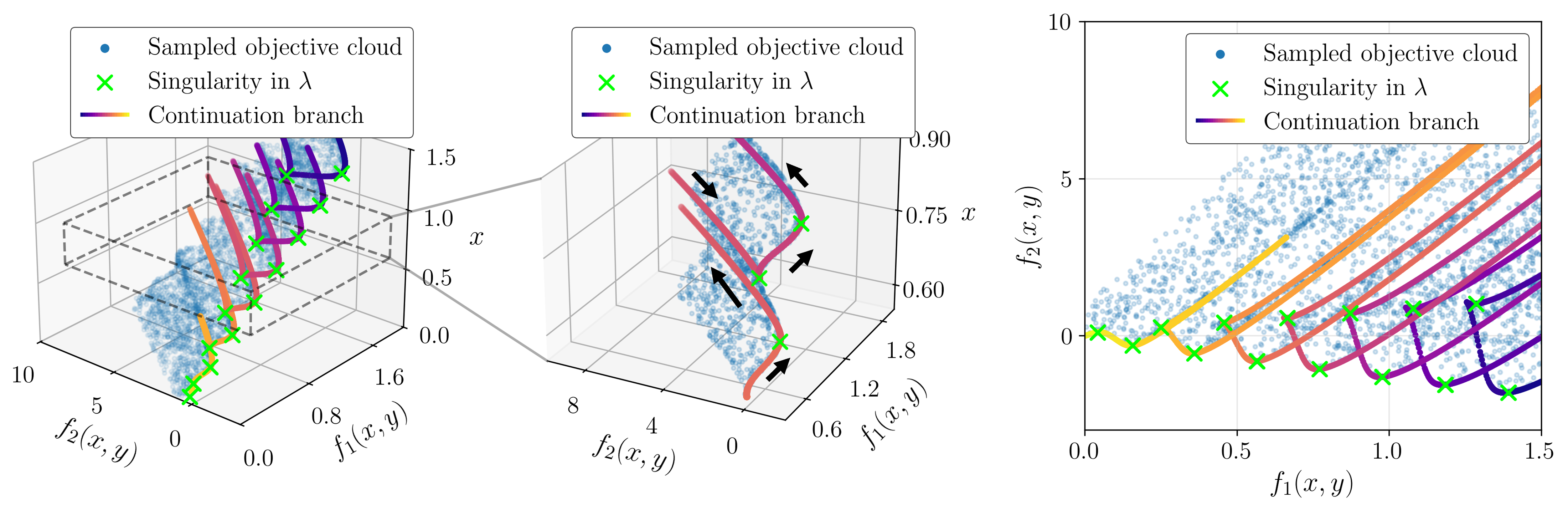}
    \caption{
    Modified two-dimensional ZDT3 with $f_1(x,y)=x+y$; continuation in $(x,\omega,\theta)$. Insets magnify the regions of large variation in $y$ (top) and the $|\lambda|$-singularities (bottom).
    The lower row visualizes the image of the continuation branch in the augmented objective space $(f_1,f_2,x)$ (left), along with a magnified, rotated view (center) that shows how the continuation trajectory (black arrows). The corresponding projection into the objective space $(f_1,f_2)$ is shown on the right.
    }
    \label{fig:zdt3_2d_modified}
\end{figure}

\begin{figure}[H]
    \centering
    \includegraphics[width=0.75\linewidth]{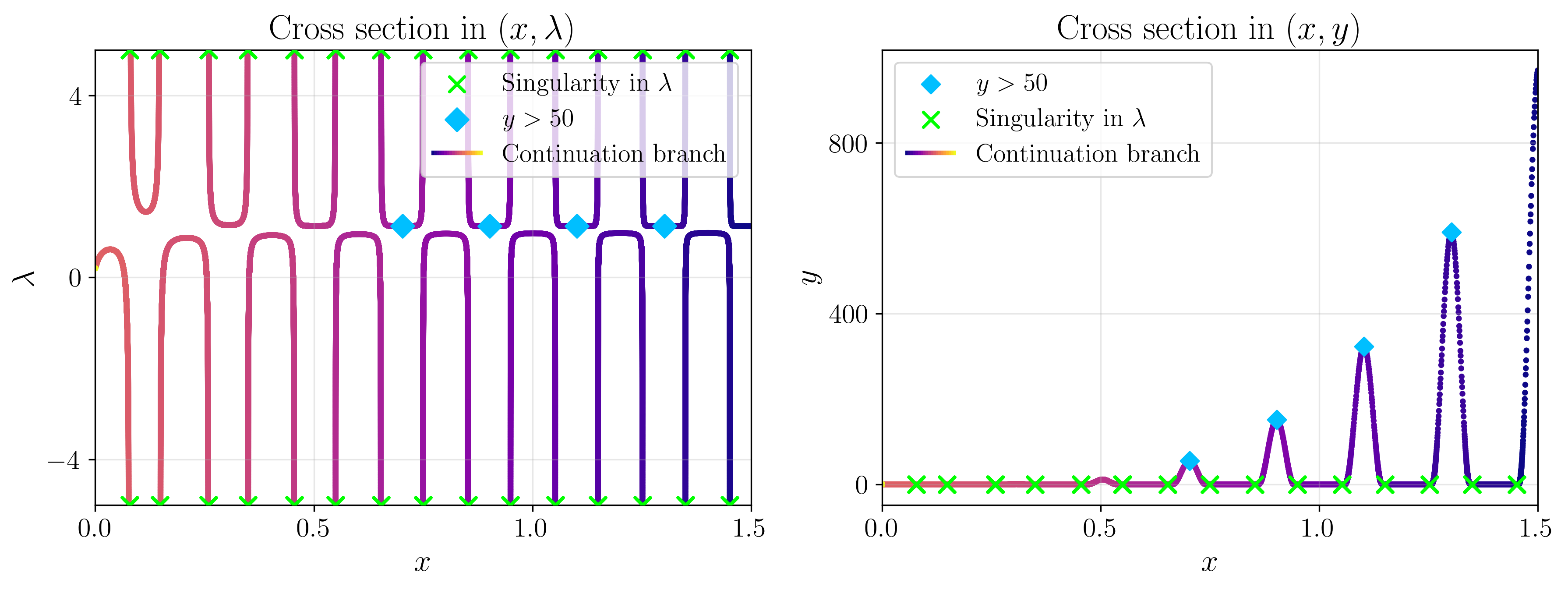}
    \caption{Cross-sections of the continuation branch for the modified two-dimensional ZDT3 with $f_1(x,y)=x+y$;. \emph{Left:} $(x,\lambda)$ showing repeated $\lambda$-singularities (green) and large values of $y$ (blue). \emph{Right:} $(x,y)$ showing the same branch. As the successive peaks in $y$ grow rapidly, the corresponding compactified coordinates approach their respective boundaries, $\omega\to\pi/2$ under angular compactification and $Y\to1$ under projective compactification.}
    \label{fig:zdt3_modified_sections}
\end{figure}

\begin{figure}[H]
    \centering
    \includegraphics[width=0.8\linewidth]{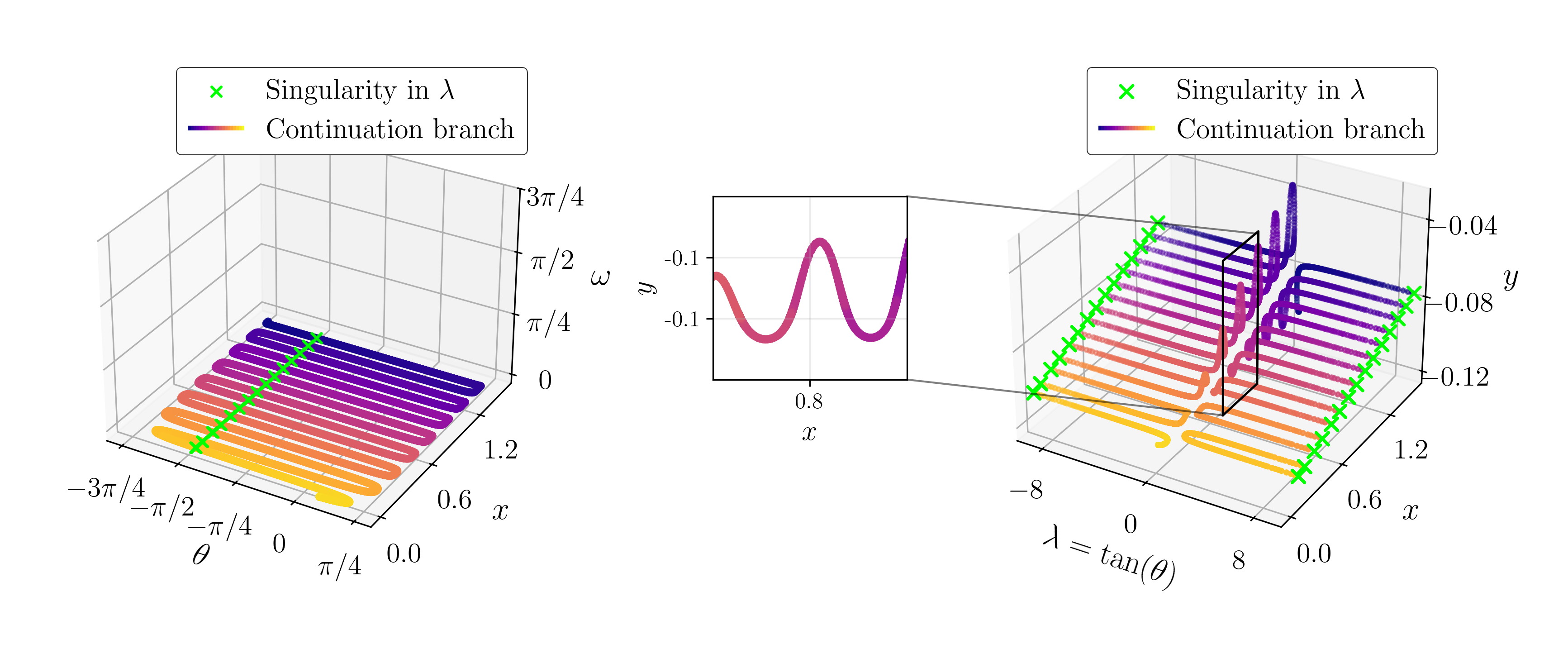}
    \includegraphics[width=0.8\linewidth]{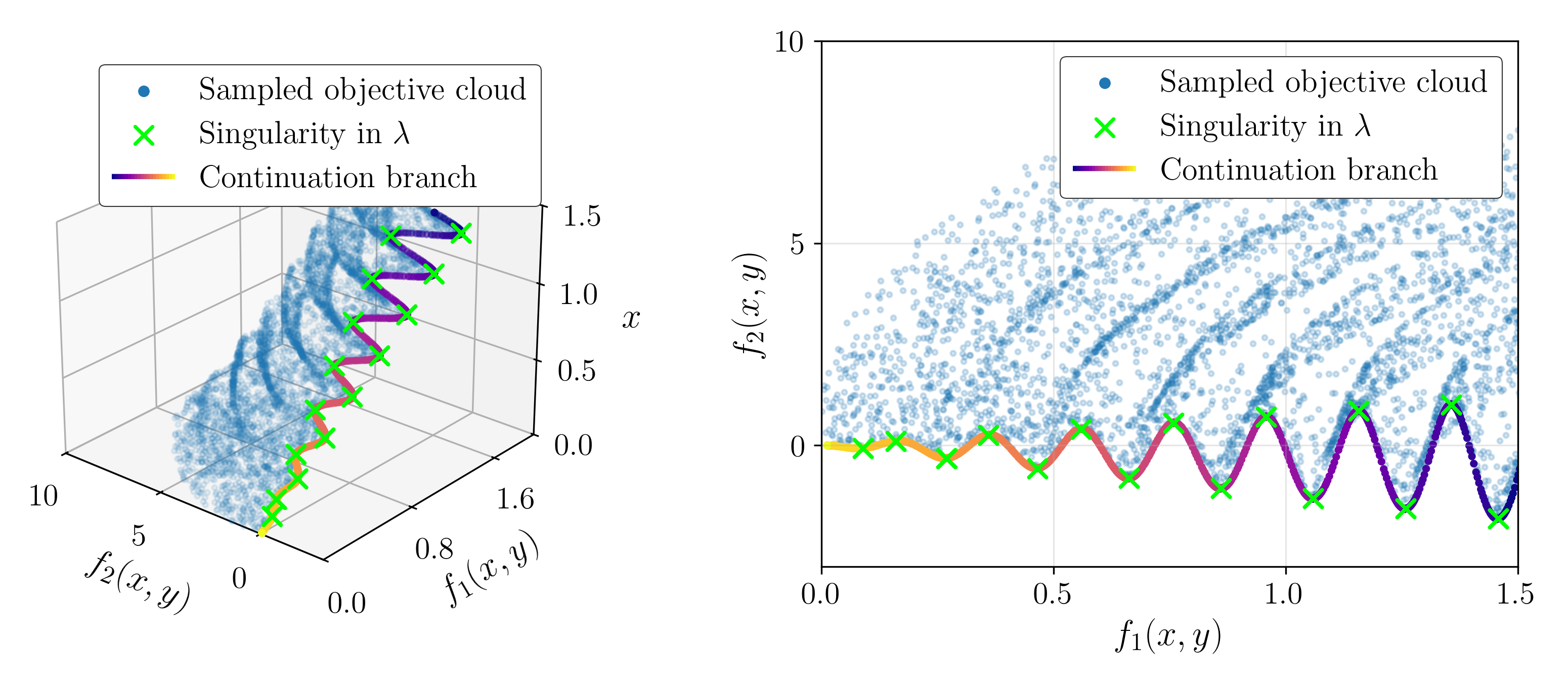}
    \caption{
    Modified two-dimensional ZDT3 with $f_1(x,y)=x+y^2$; continuation in $(x,\omega,\theta)$. Inset highlights the x-y plane near $|\lambda|$-singularities. 
    The lower row visualizes the image of the continuation branch in the augmented objective space $(f_1,f_2,x)$ (left) providing a three-dimensional view of the folds in the objective-space projection.. The corresponding projection into the objective space $(f_1,f_2)$ is shown on the right.
    }
    \label{fig:zdt3_2d_modified2}
\end{figure}

We further illustrate the angular compactification on two modified ZDT3 problems. The first modification,
\[
f_1(x,y)=x+y,
\]
was introduced and discussed in the main text. Here, we complement those results by showing the complete continuation branch in $(x,\omega,\theta)$ coordinates (Fig.~\ref{fig:zdt3_2d_modified}) together with cross-sections in $(x,\lambda)$ and $(x,y)$ (Fig.~\ref{fig:zdt3_modified_sections}), highlighting the repeated $\lambda$-singularities and the progressively larger excursions in the decision variable $y$.
We also consider a second modification,
\[
f_1(x,y)=x+y^2,
\]
while retaining the same $f_2$ and extended domain
\[
x\in[0,1.5],\qquad
y\in\left[-\tfrac19,\infty\right].
\]
Figure~\ref{fig:zdt3_2d_modified2} demonstrates that continuation in the compactified coordinates $(x,\omega,\theta)$ again traverses all apparent $\lambda$-singularities, illustrating that the angular compactification is not specific to a particular choice of objective function.

\subsection{Modified three-dimensional ZDT3}\label{app:zdt_3d_mod}
\begin{figure}[H]
    \centering
    \includegraphics[width=\linewidth]{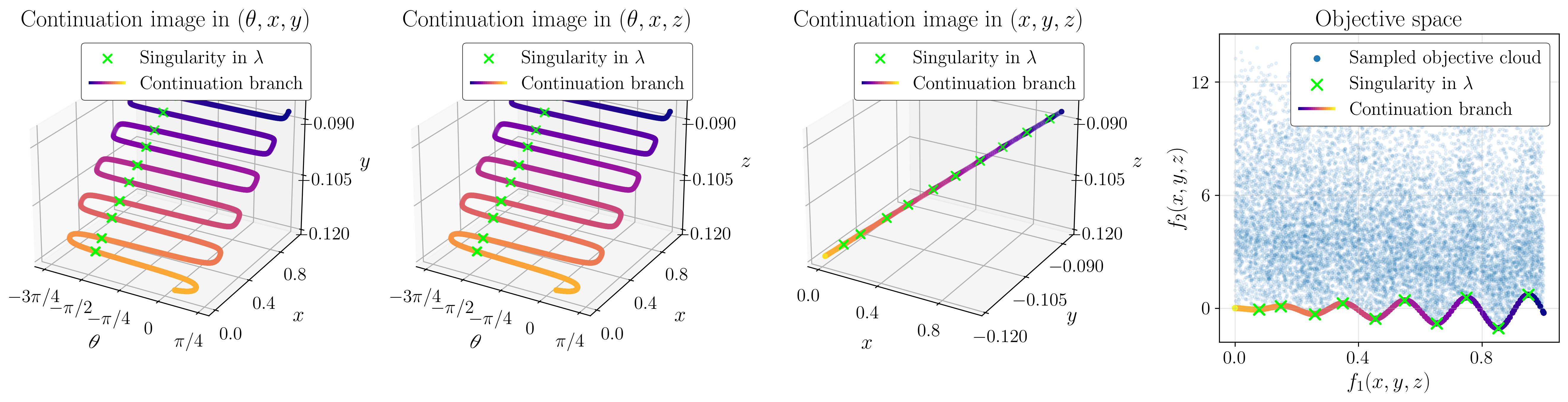}
    \caption{Modified three-dimensional ZDT3 problem. Continuation in $(x,y,z,\theta)$ with $\lambda=\tan\theta$ traverses the $|\lambda|$-singularities; the image in $(x,y,z,\lambda)$ fragments.}
    \label{fig:zdt3_3d_modified}
\end{figure}

We consider a modified 3D ZDT3 formulation with objectives
\begin{eqnarray}
    f_1(x,y,z)&=&x, \\
    f_2(x,y,z)&=&\ell(x,y,z)\left(1-\sqrt{\tfrac{x}{\ell(x,y,z)}}-\tfrac{x\sin(10\pi x)}{\ell(x,y,z)}\right),
\end{eqnarray}
where $\ell(x,y,z)=1+\tfrac{9}{2}(y+z)+5yz$. The continuation variables are $(x,y,z,\theta)$ with $\lambda=\tan\theta$, and the stationarity equations are $\nabla_{x,y,z}[\lambda f_1+(1-\lambda)f_2]=0$. Unlike the standard 3D ZDT3, the interaction term $yz$ couples the stationary equations fully, removing the rank degeneracy of the original extension and yielding a genuine 1D continuation manifold in $(x,y,z,\theta)$. Continuation on this compactified manifold again traverses the apparent singularities $|\lambda|\to\infty$ smoothly, while the image in $(x,y,z,\lambda)$ shows infinity crossings and disconnected projections (Fig.~\ref{fig:zdt3_3d_modified}).

These higher-dimensional experiments illustrate that the compactification framework remains applicable as the stationary geometry evolves from isolated curves into higher-dimensional continuation manifolds.

\section{Angular compactification for problems with $f_1'' \neq0$}\label{app:tantheta}
\begin{figure}[H]
\centering
\includegraphics[width=\linewidth]{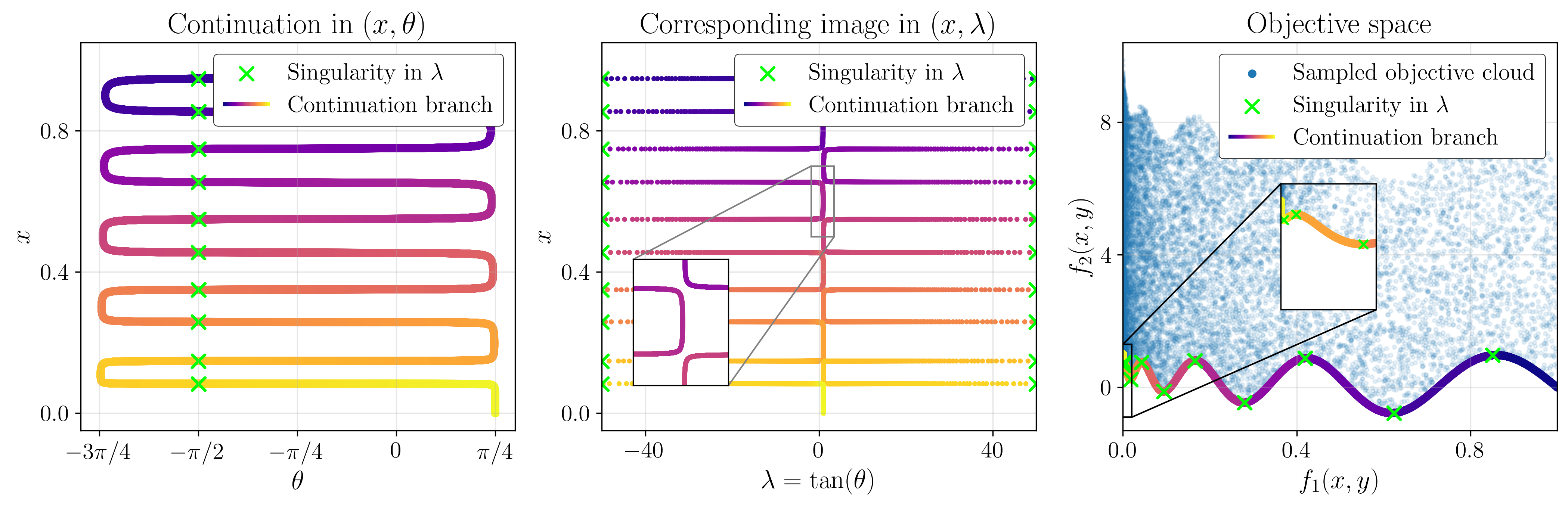}
\caption{Modified ZDT3 with $f_1(x)=x^3$.
 \emph{Left:} smooth connected branch in $(\theta,x)$. \emph{Center:} the same branch in $(\lambda,x)$ appears fragmented. The inset highlights disconnections in $(x, \lambda)$. \emph{Right:} projection into objective space. Green markers locate $|\lambda|$-singularities. The inset magnifies the images of these singularities in objective space, which become tightly clustered near $f_1=0$.}
\label{fig:zdt3x3}
\end{figure}

Unlike the standard ZDT3 problem \cite{deb2001}, for a general $f_1$ with $f_1''(x)\neq 0$ the chain rule gives
\begin{equation}
\frac{dx}{d\theta} = -\frac{f_2'(x)\,(1+\lambda^2)}{\lambda^2 f_1''(x)+\lambda(1-\lambda)f_2''(x)}
\;\xrightarrow{|\lambda|\to\infty}\; \frac{f_2'(x^\star)}{f_2''(x^\star)-f_1''(x^\star)} = -\frac{1}{c},
\end{equation}
finite, since $1+\lambda^2\sim\lambda^2$ and $1-\lambda\sim-\lambda$. Although $\sec^2\theta$ and the denominator both diverge as $\theta\to\pm\pi/2$, their ratio is bounded---an $\infty/\infty$ form resolved by leading-order cancellation. Figs.~\ref{fig:zdt3x3}--\ref{fig:zdt3sinx} verify this on modified ZDT3 with $f_1(x)=x^3$ and $f_1(x)=\sin x$, where the simplifications $f_1'=1$, $f_1''=0$ no longer hold.

\begin{figure}[H]
\centering
\includegraphics[width=\linewidth]{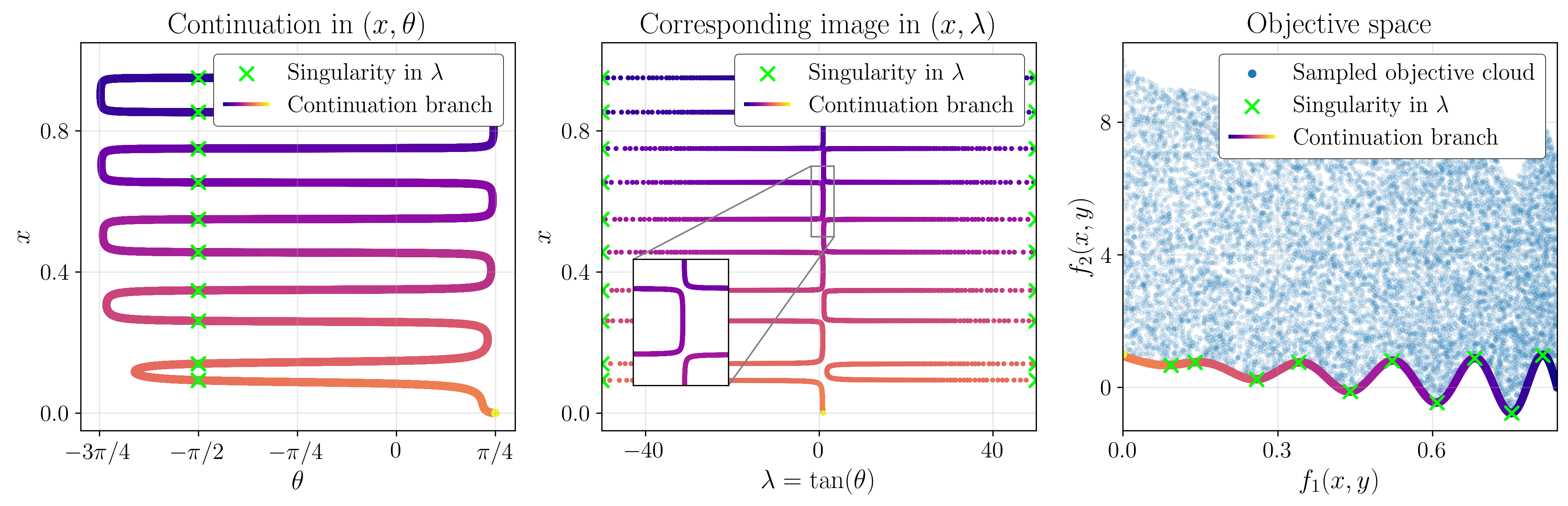}
\caption{Modified ZDT3 with $f_1(x)=\sin x$. Panels as in Fig.~\ref{fig:zdt3x3}.}
\label{fig:zdt3sinx}
\end{figure}

\section{Decision-space \underline{and} projective compactification}\label{app:decision_proj}
The compactification and reparameterization strategies presented in the main manuscript are not restricted to the scalarization parameter $\lambda$; similar singularities may arise directly in \textit{the decision variables} themselves. 
Consider the following bi-objective problem \cite{lovison2011}
\begin{equation}\label{eq:appbiobj}
f_1(x,y)=x^2+y^2, \qquad f_2(x,y)=(x-6)^2-(y+0.3)^2.
\end{equation}
Here, continuation can produce singluar behavior in the decision variable $y$, analogous to the divergence observed in $\lambda$. The corresponding stationary branch is:
\begin{equation}\label{eqn:x_stationary}
x = 6(1-\lambda),
\qquad
y = \frac{0.6(1-\lambda)}{4\lambda-2},   
\end{equation}
which becomes singular at $\lambda=\frac12$ as $y\to \pm \infty$.
To handle this, we introduce a compactified angular coordinate:
\begin{equation}
    y=\tan(\omega),
\end{equation}
and perform continuation in $(x,\omega,\lambda)$ rather than $(x,y,\lambda)$. This maps the unbounded variable $y$ onto the compact coordinate $\omega$, so that the singular points $|y|\to\infty$ correspond to the finite values $\omega=\pm\frac{\pi}{2}$.

As with the $\lambda=\tan(\theta)$ transformation, the resulting compactified manifold remains smooth and traversable across the apparent singularity; branches that appear disconnected in $(x,y,\lambda)$ become connected components of a single smooth manifold in the compactified coordinates.
This suggests that many singularities encountered in continuation-based multi-objective optimization are coordinate artifacts rather than intrinsic geometric features, and that compactification provides a systematic way to regularize them.

\begin{figure}[H]
    \centering
    \includegraphics[width=\linewidth]{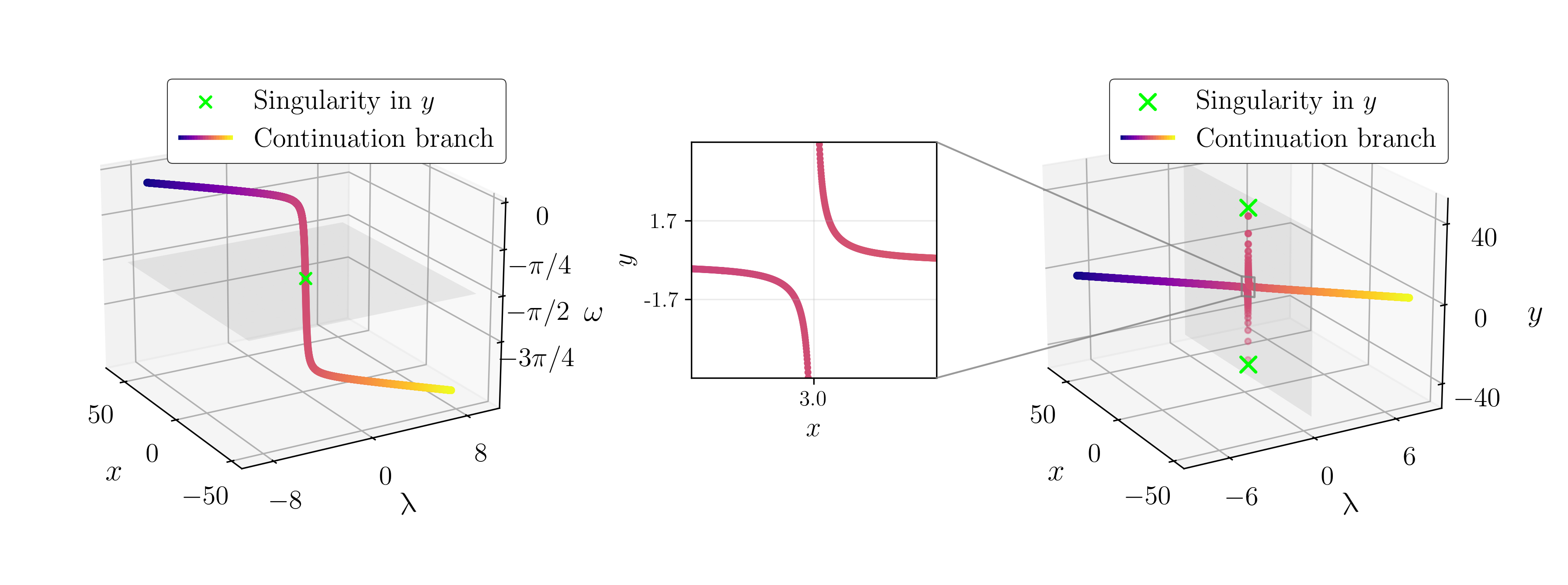}
    \includegraphics[width=\linewidth]{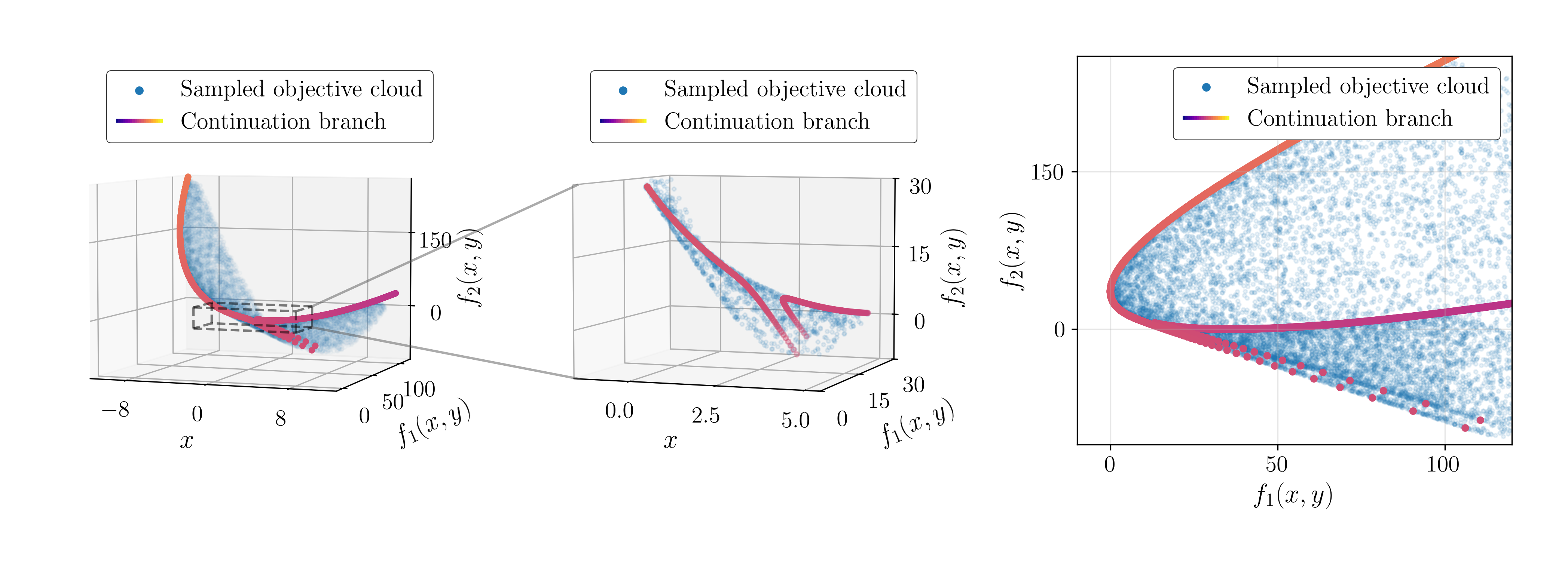}
    \caption{
    Continuation through a decision-space singularity using angular compactification. The continuation branch is shown in the compactified coordinates $(x,\omega,\lambda)$ (top left) and mapped back to the original variables $(x,y,\lambda)$ (top right). The insets magnifies the singularity in $(x,\omega)$. The transformation $y=\tan\omega$ regularizes the singularity at $y=\pm\infty$ by mapping it to the finite point $\omega=\pm\pi/2$, enabling pseudo-arclength continuation to traverse the decision-space infinity. The lower row visualizes the image of the continuation branch in the augmented objective space $(x,f_1,f_2)$ (left), where the additional coordinate $x$ reveals the folded geometry obscured by the objective-space projection. A magnified, rotated view (center) highlights how branches that appear superimposed in objective space separate. The corresponding projection into the objective space $(f_1,f_2)$ is shown on the right.
    }
    \label{fig:x_singular}
\end{figure}

\subsection{Projective Compactification: Single-scale compactification}
Consider the case where one or more variables may become unbounded along the branch. When the singular behavior is associated with a single point at infinity, the entire state can be compactified using a common projective denominator. For a two-dimensional state $\bm{x}=(x,y)$ and continuation parameter $\lambda$, we introduce the coordinates
\begin{equation}
    X = \rho x
\qquad
    Y = \rho y,
\qquad
L = \rho \lambda,
\end{equation}
together with the normalization constraint
\begin{equation}
    X^2+Y^2+L^2+\rho^2=1. 
    %\qquad \rho = \frac{1}{\sqrt{1 + x^2 + y^2 + \lambda^2}}    
\end{equation}
This maps the unbounded space $(x,y,\lambda)\in\mathbb{R}^3$ onto the compact unit sphere in $(X,Y,L,\rho)$ coordinates. 
%The physical variables are recovered through
%x=\frac{X}{\rho}, \qquad y=\frac{Y}{\rho}, \qquad \lambda=\frac{L {\rho}.
Points satisfying $\rho=0$ correspond to points at infinity in the original coordinates. Unlike the physical representation, these points remain finite in the compactified space and can therefore be traversed directly by continuation algorithms.
As an illustrative example, consider the unconstrained bi-objective problem in Eq.~(\ref{eq:appbiobj}). Substituting the projective coordinates in Eq.~(\ref{eqn:x_stationary}) and multiplying through by appropriate powers of $\rho$ yields the homogenized system
\[
X-6(\rho-L)=0,
\]
\[
(2L-\rho)Y-0.3(\rho-L)\rho=0.
\]
Pseudo-arclength continuation can therefore proceed smoothly across $\rho=0$ which, in the original coordinates, corresponds to traversing infinity and connecting asymptotic branches that would otherwise appear separated at $+\infty$ and $-\infty$. For problems in which all singular variables share a common asymptotic scaling and a unique singular location, this shared-denominator formulation provides a simple and effective regularization strategy.

\begin{figure}[H]
    \centering
    \includegraphics[width=\linewidth]{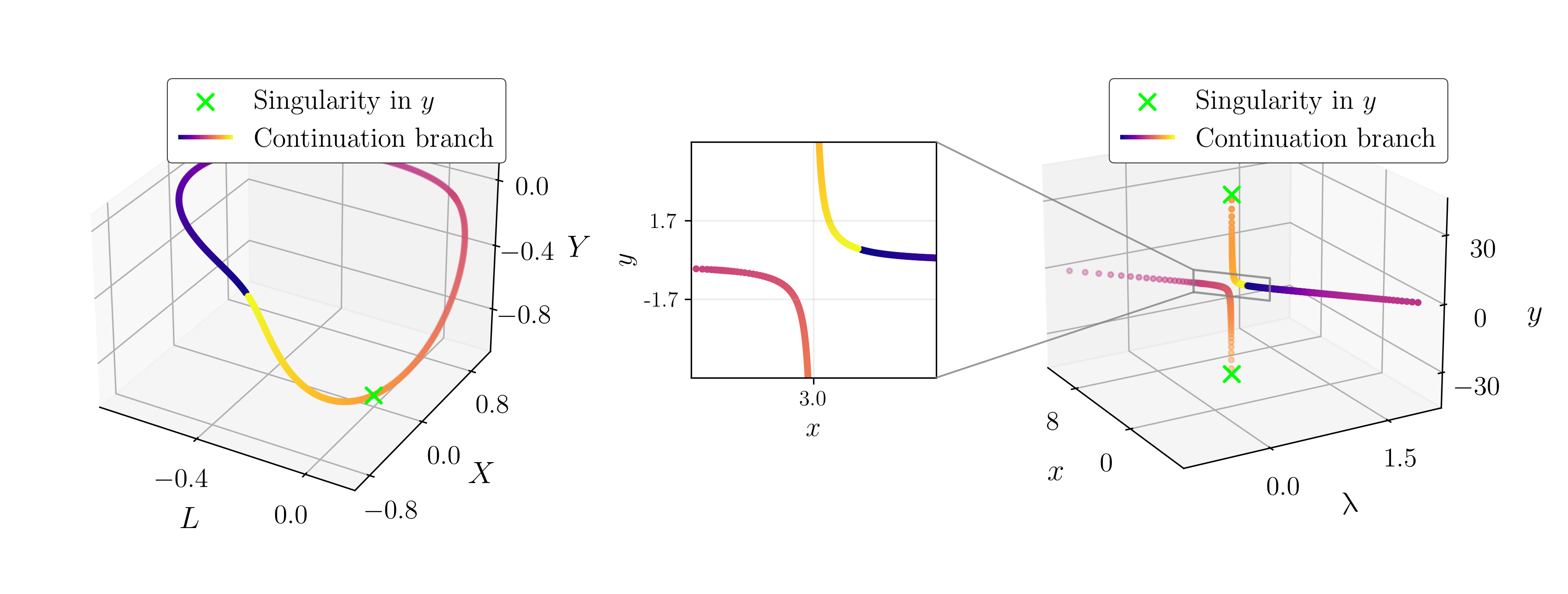}
    \includegraphics[width=\linewidth]{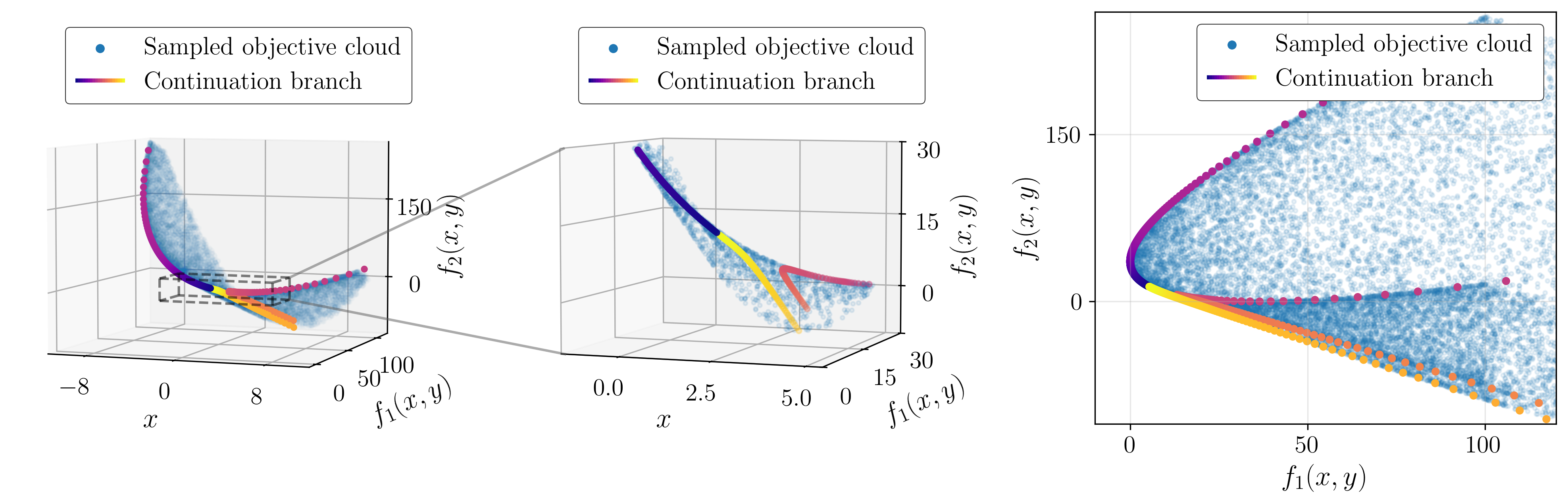}
    \caption{
    Continuation through the decision-space singularity using single-scale projective compactification, along the coordinates $(X,Y,L)$ with shared denominator $\rho$. The inset highlights the disconnection in $(x, \lambda)$. The lower row visualizes the image of the continuation branch in the augmented objective space $(x,f_1,f_2)$. Panels as in Fig.~\ref{fig:x_singular}.
    }
    \label{fig:x_singular_c}
\end{figure}
\end{document}